\documentclass[11pt,reqno,a4paper,twoside]{article}
\usepackage{amsmath,amsthm,amstext,amscd,amssymb,euscript}
\renewcommand{\phi}{\varphi}

\renewcommand{\subset}{\subseteq}
\renewcommand{\emptyset}{\varnothing}

\def\limn{\lim_{n\to\infty}}
\newcommand{\Zd}{\mathbb Z^d}

\renewcommand{\Pr}{\mathbb P}

\def\1{ {\mathit{1} \!\!\>\!\! I} }

\makeatletter
\@addtoreset{equation}{section}
\makeatother

\newtheorem{ittheorem}{Theorem}
\newtheorem{itlemma}{Lemma}
\newtheorem{itproposition}{Proposition}
\newtheorem{itdefinition}{Definition}
\newtheorem{itremark}{Remark}

\newenvironment{theorem}{\addtocounter{equation}{1}
\begin{ittheorem}}{\end{ittheorem}}

\newenvironment{lemma}{\addtocounter{equation}{1}
\begin{itlemma}}{\end{itlemma}}

\newenvironment{proposition}{\addtocounter{equation}{1}
\begin{itproposition}}{\end{itproposition}}

\newenvironment{definition}{\addtocounter{equation}{1}
\begin{itdefinition}}{\end{itdefinition}}

\newenvironment{remark}{\addtocounter{equation}{1}
\begin{itremark}}{\end{itremark}}

\newcommand{\beq}{\begin{eqnarray}}
\newcommand{\eeq}{\end{eqnarray}}

\newcommand{\be}{\begin{equation}}
\newcommand{\ee}{\end{equation}}

\newcommand{\bl}{\begin{lemma}}
\newcommand{\el}{\end{lemma}}

\newcommand{\br}{\begin{remark}}
\newcommand{\er}{\end{remark}}

\newcommand{\bt}{\begin{theorem}}
\newcommand{\et}{\end{theorem}}

\newcommand{\bd}{\begin{definition}}
\newcommand{\ed}{\end{definition}}

\newcommand{\bp}{\begin{proposition}}
\newcommand{\ep}{\end{proposition}}

\newcommand{\bc}{\begin{corollary}}
\newcommand{\ec}{\end{corollary}}

\newcommand{\bpr}{\begin{proof}}
\newcommand{\epr}{\end{proof}}

\newcommand{\bi}{\begin{itemize}}
\newcommand{\ei}{\end{itemize}}

\newcommand{\ben}{\begin{enumerate}}
\newcommand{\een}{\end{enumerate}}

\newcommand{\Z}{\mathbb Z}
\newcommand{\R}{\mathbb R}
\newcommand{\N}{\mathbb N}
\newcommand{\Comp}{\mathbb C}
\newcommand{\C}{\mathcal C}
\newcommand{\Q}{\mathbb Q}
\newcommand{\E}{\mathbb E}

\newcommand{\bee}{\ensuremath{\mathcal{B}}}
\newcommand{\gee}{\ensuremath{\mathcal{G}}}

\newcommand{\s}{\ensuremath{\mathcal{S}}}

\newcommand{\fe}{\ensuremath{\mathcal{F}}}

\newcommand{\si}{\ensuremath{\sigma}}
\newcommand{\epsi}{\ensuremath{\epsilon}}

\def\now{
\ifnum\time<60       
          12:\ifnum\time<10 0\fi\number\time am 
          \else
            \ifnum\time>719\chardef\a=`p\else\chardef\a=`a\fi 
          \hour=\time
          \minute=\time 
          \divide\hour by 60 
          \ifnum\hour>12\advance\hour by -12\advance\minute by-720 \fi
          \number\hour:%
          \multiply\hour by 60 
          \advance\minute by -\hour
          \ifnum\minute<10 0\fi\number\minute\a m\fi}           
\newcount\hour 
\newcount\minute
\numberwithin{equation}{section}         

\theoremstyle{remark}

\def\t{{\bf t}}  
\def\r{{\bf r}}  
\def\w{{\bf w}}  

\begin{document}

\title{{\bf Exponential distribution
for the occurrence of rare patterns
in Gibbsian random fields }}

\author{
M. Abadi
\footnote{IME-USP, cp 66281, 05315-970, S{\~a}o Paulo, SP, Brasil,
abadi@ime.usp.br.}\\
J.-R. Chazottes
\footnote{CPhT, CNRS-Ecole polytechnique, 91128 Palaiseau Cedex, France,
jeanrene@cpht.polytechnique.fr}\\
F.\ Redig
\footnote{Faculteit Wiskunde en Informatica, Technische Universiteit
Eindhoven, Postbus 513,
5600 MB Eindhoven, The Netherlands, f.h.j.redig@tue.nl}\\
E. Verbitskiy
\footnote{Philips Research Laboratories, Prof. Holstlaan 4, 5656 AA
Eindhoven, The Netherlands,
evgeny.verbitskiy@philips.com}
}

\maketitle

\footnotesize
\begin{quote}
{\bf Abstract:}
We study the distribution of the occurrence of rare patterns
in sufficiently mixing Gibbs random fields on the lattice
$\mathbb{Z}^d$, $d\geq 2$.
A typical example is the high temperature Ising model.
This distribution is shown to converge to an exponential law as the
size of the pattern diverges. Our analysis not only provides this
convergence
but also  establishes a precise estimate of the distance between the exponential law
and the distribution of the occurrence of finite patterns.
A similar result holds for the
repetition of a rare pattern. We apply these results to the fluctuation
properties of occurrence and repetition of patterns:
We prove a central limit theorem and a large deviation principle.
\end{quote}
\normalsize


\vspace{12pt}
\section{Introduction}
In the last decade there has been an intensive study of exponential laws
for rare events
in the context of dynamical systems and stochastic processes, see e.g. the
review
paper \cite{AG}. In general, these laws are derived under the assumption
of sufficiently
strong mixing conditions, which basically ensures the possibility of
writing the rare
event as an intersection of almost independent events. The basic example
of a rare event is the
occurrence or return of a large cylindrical event. Other relevant examples
are approximate
cylindrical events (approximate matching in the sense of Hamming
distance, see e.g. \cite{Chi}), or large deviation
events in certain interacting particle systems, see e.g.
\cite{Asselah,Asselahbis}.

The mixing conditions appearing in the context of dynamical systems or
stochastic processes are typical
for $\Z$-actions, e.g., the $\psi$-mixing condition is very naturally
satisfied in the context
of Bowen-Gibbs measures \cite{bowen}. In turning to the context of random
fields
or $\Z^d$-actions, the $\psi$-mixing property is very restrictive and in
many natural examples
such as Gibbsian random fields, this property does not hold except
(trivially) in the i.i.d. case and in non-interacting
copies of one-dimensional Gibbs measures.

Gibbsian random fields have an obvious relevance to various
applications, e.g., statistical physics, image processing,
etc.
Many interesting fluctuation properties such as large deviations
principle,
central limit theorems have been derived for them, and
by now Gibbs measures constitute
a well-established field of research, see e.g. \cite{ellis}, \cite{Geo},
\cite{Guyon}. 

The study of exponential laws for the occurrence or repetition of rare
events in random fields
has been initiated by A.J. Wyner \cite{wyner} for the $\psi$-mixing
case, using the Chen-Stein method. Because
of the mixing condition, 
the results of that paper are not applicable
to Gibbsian random fields like the Ising model in the high mixing
regime (such as Dobrushin uniqueness, or analyticity regime).

As an example, consider the $d$-dimensional Ising model in the
high-temperature regime
and fix a pattern in a cubic box of size $n$: what is the size of the
``observation window" in which we
see this pattern for the first time ? This is clearly a rare event when
the size of the pattern increases, and
hence one expects in the ``high mixing regime"
that the size of this observation window is approximately exponentially
distributed with parameter
proportional to the probability of the pattern.

The main difficulty in making this intuition into a mathematical
statement is caused by the typical non-uniform
mixing of Gibbsian random fields: the influence of an event $A$ on
an event $B$  is not only dependent
on their distance but also on their size. More precisely, the difference
between the conditional
probabilities $\Pr (A|B)$ and $\Pr (A)$ can be estimated in the optimal
situation of Dobrushin uniqueness regime
as something of the form $|A|\exp(-\textup{dist}(A,B))$.
On a technical level, this ``non-uniform mixing" implies that the rare
event
under consideration should be written as an intersection of events which
at the same time
are separated by a large distance and do not have an ``excessive" size. 

In this paper we concentrate on Gibbsian random fields in the Dobrushin
uniqueness regime (e.g. high temperature
case). This has to be considered as the first non-trivial test case for
random fields, with a broad
variety of examples. The regime of phase coexistence (such as in the
low-temperature Ising model) poses
an even larger non-uniformity in the mixing conditions, i.e., the
difference between $\Pr (A|B)$ and
$\Pr (A)$ will in that case also depend on which events $B$ we are
conditioning on. Recent techniques
such as disagreement percolation constitute a powerful tool to tackle this
situation. 
This is however not the subject of the present paper, where we want to
deal with the basic non-uniformity in
the mixing appearing in all non-trivial Gibbsian random fields.

Besides the mere derivation of exponential laws for the occurrence and
repetition of rare events,
we obtain a precise and uniform 
estimate of the error (i.e., the difference between the
law and its exponential approximation). We show that obtaining this
precise control of the error has many
useful non-trivial applications in studying fluctuations of both "waiting
times" and repetitions of rare patterns.
The derivation of the exponential
law is not via Chen-Stein method. Via a direct
use of the (non-uniform) mixing
we obtain more detailed information on the error term. The reason for that
is that in the  Chen-Stein method one  gives an estimate
of the variational distance between
the ``real counting process" and the Poisson process, 
whereas we only need one particular event. The precise
estimation of the error turns to be crucial in the study of large deviations.

The problem of ``waiting times" is to ask for
the $\Pr$-typical size of the ``observation window" in which a
$\Q$-typical pattern occurs, where
$\Pr$ is Gibbsian, and $\Q$ is any ergodic field. The logarithm of the
size of this observation
window properly normalized converges to the sum of the entropy of $\Q$ and
the relative entropy
density $s(\Q|\Pr)$. To this ``law of large numbers" we add precise large
deviation estimates and a central limit
theorem as a corollary
of the exponential law with its precise error. The main point is that the
exponential law provides an
approximation of the logarithm of
the waiting time by minus the logarithm of the probability of the
corresponding pattern.
For the cumulant generating function of the
waiting times, we give an explicit expression in terms of the
pressure. It coincides with the cumulant generating function of
the probability of patterns in the interval $(-1,\infty)$ and
is constant on $(-\infty, -1]$.
A similar phenomenon was observed numerically for the cumulant generating
function of the return times (that is in dimension
one), see \cite{vaienti}.

For repetition of patterns, we prove a similar exponential law
with precise error bound. However, in that case we have to exclude
``badly self-repeating" patterns, which have exponentially small
probability for any Gibbs measure. As a corollary, we obtain a
law of large numbers and a central limit theorem for
repetitions. The large deviations are more subtle due to the
presence of the bad patterns. We prove a full large deviation
principle for the measure conditioned on good patterns, and
a restricted large deviation principle for the full measure.

Our paper is organized as follows. In section 2 we give basic notations
and definitions and state
our main result and its corollaries. In section 3 we review basic
properties of high-temperature Gibbs measures. 
Section 4 contains the proof of the exponential law for the occurrence of
patterns, and section 5
is devoted to the derivation of its corollaries.

\section{Definitions and results}

We consider a random field $\{\sigma({\bf x}):{\bf x}\in\Z^d\}$
on the lattice $\Z^d$, $d\geq 2$, where $\sigma({\bf x})$ takes values
in a finite set $\mathcal{A}$. The joint distribution
of $\{ \sigma({\bf x}):{\bf x}\in\Z^d\}$ is denoted
by $\Pr$.
The configuration space
$\Omega = \mathcal{A}^{\Z^d}$ is endowed with the product
topology (making it into a compact metric space).
The set of finite subsets of $\Z^d$ is denoted by $\s$.
For $A,B\in\s$ we put $d(A,B)=\min\{ |{\bf x}-{\bf y}|:{\bf x}\in A,{\bf
y}\in B\}$, where
$|{\bf x}|=\sum_{i=1}^d |x_i|$ (${\bf x}=(x_1,x_2,...,x_d)$).
For $A\in\s$, $\fe_A$ is the sigma-field 
generated by $\{\sigma({\bf x}):x\in A\}$.
For $V\in\s$ we put $\Omega_V=\mathcal{A}^V$. For $\sigma\in\Omega$,
and $V\in\s$, $\si_V\in\Omega_V$ denotes the restriction of
$\si $ to $V$. For ${\bf x}\in\Z^d$ and $\si\in\Omega$,
$\tau_{{\bf x}} \si$ denotes the translation of $\si$ by ${\bf x}$:
$\tau_{{\bf x}}\si ({\bf y}) = \si({\bf x}+{\bf y})$.
For an event $E\subset\Omega$ the dependence
set of $E$ is the minimal $A\in\s$ such that
$E$ is $\fe_A$ measurable.
For any $n\in\N$ let $C_n= [0,n]^d\cap\Z^d$.
An element $A_n\in \Omega_{C_n}$ is called a {\bf $n$-pattern} or
a pattern of size $n$.

\bd[First occurrence of a pattern]
For every configuration $\sigma\in\Omega$ 
we define $\t_{A_n}(\sigma)$ to be the first occurrence of an $n$-pattern
$A_n$ in that configuration,
that is the minimal $k\in \N$ such that there exists a non-negative vector 
${\bf x}= (x_1,\ldots,x_d)\in\Z_{+}^d$ with $x_i\le k$, $i=1,...,d$,
$|{\bf x}|>0$, satisfying
\be
(\tau_{{\bf x}} \sigma)_{C_n}=  A_n.
\ee
If such a vector $x$ does not exist then we put $\t_{A_n} (\si ) =\infty$.
\ed

We now come to the mixing hypothesis we make on our random fields. For
$m>0$ define 
\be\label{mixfun}
 \phi(m) = \sup \frac 1{|A_1|}\, |\ \Pr\left(E_{A_1}|E_{A_2}\right) -
 \Pr\left(E_{A_1}\right)|\, ,
\ee
where the supremum is taken over all finite subsets $A_1, A_2$ of $\Zd$,
with $d(A_1,A_2)\ge m$ and 
$E_{A_i}\in\fe_{A_i}$, with $\Pr(E_{A_2})>0$.
Note that this $\phi (m)$ differs from the usual
$\varphi$-mixing function since we divide by the size of the dependence
set of the event $E_{A_1}$.
\bd\label{mixdef} 
A random field is  {\bf non-uniformly exponentially $\phi$-mixing}
if there exist  constants $C_1, C_2 >0$ such that 
\be\label{exponmix}
  \phi(m) \le C_1 e^{-C_2 m}\quad \textup{for all}\quad m>0.
\ee
\ed

The examples that motivate this definition are
Gibbsian random fields in the Dobrushin uniqueness regime (see
Definition \ref{dobrushin} below and examples thereafter).
We leave their definition and properties till the next section.
For a pattern $A_n\in \Omega_{C_n}$ we define the corresponding cylinder
$\C(A_n)$ as 
$$
\C(A_n)=\{\sigma\in\Omega: \sigma_{C_n}=A_n\}\,.
$$


Our main result reads:

\bt\label{thm1}
For a translation-invariant Gibbs random field satisfying
\eqref{exponmix},
there exist strictly positive constants $C,c,\rho,\Lambda_1,\Lambda_2$,
$\Lambda_1\leq\Lambda_2$,
such that for any $n$ 
and any $n$-pattern $A_n$, there exists
$\lambda_{A_n}\in[\Lambda_1,\Lambda_2]$, 
such that
\be\label{main}
\Bigl|\ \Pr\ \Bigl\{
\t_{A_n}>\Bigl(\frac{t}{\lambda_{A_n}\Pr\left(\mathcal{C}(A_n)\right)}\Bigr)^{1/d}\Bigr\}
-e^{-t}\ \Bigr| \le\ C\ \Pr\left(\mathcal{C}(A_n)\right)^\rho\ e^{-ct}
\ee
for any $t>0$.
\et
Notice that $\Pr(\mathcal{C}(A_n))$ in the ``error term'' in (\ref{main})
is bounded above
by $\exp(-c' n^d)$, with $c'>0$, 
by the Gibbs property, see (\ref{gibbsineq}).

The proof of this theorem is given in Section \ref{mainproof}.

\br The only results we are aware of in the context of random fields
appeared in \cite{wyner}.
The results of that paper are valid under the assumption of a much
stronger mixing condition
than ours, namely $\psi$-mixing. Most Gibbs random fields (including
the Ising model at high temperature) cannot satisfy such a property.
As an examples of $\psi$-mixing Gibbsian random fields
(in the sense of Wyner) on
$\mathbb{Z}^2$, one can consider independent
copies of a one-dimensional Markov chain, this
gives a two-dimensional Gibbsian random field, but
without interaction in the $y$-direction.

From the technical point of view, Wyner uses the Chen-Stein method.
This leads to an estimate which for fixed pattern size
does not converge to zero as $t\to\infty$.
Here we use a different approach allowing us to get a control in
$t$ in \eqref{main}. This feature will turn to be fundamental when
we prove large deviations for waiting times, see below.
\er

From the proof of Theorem \ref{thm1} it will be clear that we can
generalize it to $(A_n)_n$'s that are finite patterns supported on a van
Hove sequence
of subsets of $\Z^d$. 

We will show elsewhere how to prove an analog of Theorem \ref{thm1}
in order to obtain the same kind of result for the
low temperature ``plus phase'' of the Ising model, where
the mixing condition of Definition \ref{mixdef} is no longer satisfied.

We now state a number of corollaries of the previous theorem.
We first consider the repetition of patterns.

\bd[First repetition of the initial pattern]
For every configuration $\sigma\in\Omega$ and for all $n\in\N$, 
we define the first repetition, denoted by $\r_n(\sigma)$,
as the minimal $k\in \N$ such that there exist a  vector 
${\bf x}= (x_1,\ldots,x_d)\in\Z_{+}^d$, with $0\leq x_i\le k$ and
 $|{\bf x}|>0$, satisfying
\be
(\tau_{{\bf x}} \sigma)_{C_n}= \sigma_{C_n}\ .
\ee
\ed

To obtain a similar result for the repetition times we
have to exclude certain patterns with ``too quick repetitions''.
We will make this notion precise later.
The following result is established in Subsection \ref{proof-coro-repet}.

\bt\label{corollary-repetition}
For a translation-invariant Gibbs random field satisfying
\eqref{exponmix},
there exist
\begin{itemize}
\item[(i)] a set $G_n$, which is a union of cylinders;
\item[(ii)]strictly positive constants $B,b, C,c,\rho$
\end{itemize}
such that for any $n\geq 1$ 
\be\label{probaofbadpatterns}
\Pr( G_n^c) \le Be^{-bn^d},
\ee
and for each $A_n$ with $\mathcal{C}(A_n) \subset G_n$
\be\label{exponentiallawforrepetitions}
\Bigl|\ \Pr\ \Bigl\{
\r_n >\Bigl(\frac {t}{\lambda_{A_n}\Pr(\C(A_n))}\Bigr)^{1/d}\Big|\
\C(A_n)\Bigr\}- e^{-t}\ \Bigr|
\le \ C\ \Pr(\C(A_n))^\rho\ e^{-ct}
\ee
for all $t>0$ and where $\lambda_{A_n}$ is given in Theorem \ref{thm1}. 
\et

Notice that the constants appearing in the previous Theorems may be
different. Nevertheless we used the same notations for the sake of
simplicity.

We denote by $s(\Pr)$ the entropy of $\Pr$ (see the next section for
the definition). The next result (proved in subsection
\ref{proof-corollary-OW}) shows how the repetition of
typical patterns allows to compute the entropy using a single ``typical''
configuration.

\bt\label{corollary-OW}
For a translation-invariant Gibbs random field satisfying
\eqref{exponmix},
there exists $\epsilon_0 >0$ such that
for all $\epsilon>\epsilon_0$
\be\label{strong-approximation}
-\epsilon\log n \leq\log\left[(\r_n(\sigma))^d\ \Pr({\mathcal
C}(\sigma_{C_n}))\right]\leq \log\log n^{\epsilon}
\quad\textup{eventually}\; \Pr\!-\!\textup{almost surely}.
\ee
In particular,
\be\label{NUEMOW}
\lim_{n\to\infty}\frac{d}{n^d}\log\r_n(\sigma)=s(\Pr)\quad\Pr-\textup{a.s}\,
.
\ee
\et

Note that \eqref{NUEMOW} is a particular case of the result by Ornstein
and Weiss in \cite{OW} 
where $\Pr$ is only assumed to be ergodic. Under our assumptions, we get
the more precise
result \eqref{strong-approximation}.

We now consider the occurrence of an  $n$-pattern drawn from some
ergodic random field in the configuration drawn
from   a possibly different  Gibbsian random field.
This is the natural $d$-dimensional analog of the waiting-time
\cite{shields}, \cite{wyner}.

\bd[``Waiting time'']
For all configurations $\xi,\sigma\in\Omega$ and for all $n\in\N$, 
we define the ``waiting time'', denoted by $\w_n(\xi,\sigma)$,
as the minimal $k\in \N$ such that there exist a non-negative vector 
${\bf x}= (x_1,\ldots,x_d)\in\Z_{+}^d$, with $0\leq x_i\le k$ and $|{\bf
x}|>0$, satisfying
\be
(\tau_{{\bf x}} \sigma)_{C_n}= \xi_{C_n}\ . 
\ee
\ed
Notice that $\w_n(\xi,\sigma)=\t_{\xi_{C_n}}(\sigma)$. We are going to
consider the situation when
$\xi$ is ''randomly chosen'' according to an ergodic random field
$\mathbb{Q}$ and $\sigma$ is
''randomly chosen'' according to a non-uniformly exponentially
$\phi$-mixing Gibbs random field $\Pr$, i.e. $(\xi,\sigma)$ is drawn with
respect
to the product measure $\mathbb{Q}\times\Pr$. We denote by
$s(\mathbb{Q}|\Pr)$ the relative entropy
of ${\mathbb Q}$ with respect to $\Pr$;
see section \ref{Gibbs} for the definition and a more explicit form.
We have the following result (proved in Subsection \ref{wtproof}):

\bt\label{waiting-time}
For a translation-invariant Gibbs random field $\Pr$ satisfying \eqref{exponmix},
and an ergodic random field $\mathbb{Q}$,
there exists $\epsilon_0 >0$ such that
for all $\epsilon>\epsilon_0$
\be
-\epsilon\log n \leq\log\left[(\w_n(\xi,\sigma))^d\ \Pr({\mathcal
C}(\xi_{C_n}))\right]\leq \log\log n^\epsilon
\ee
for $\mathbb{Q}\times\Pr$-eventually almost every $(\xi,\sigma)$.
In particular
\be\label{ASW}
\lim_{n\to\infty}\frac{d}{n^d}\log
\w_n(\xi,\sigma)=s(\mathbb{Q})+s(\mathbb{Q}|\Pr)
\quad\mathbb{Q}\times\Pr-\textup{a.s}\, .
\ee
\et
Statement \eqref{ASW} is the d-dimensional generalization of a result
obtained in \cite{chazottes}
in the case of Bowen-Gibbs measures. 
Using Theorem \ref{corollary-OW} we can rewrite \eqref{ASW}, for a
``typical'' pair $(\xi,\sigma)$, as follows:
$$
\w_n(\xi,\sigma)\approx \r_n(\xi) \exp((n^{d}/d) s(\mathbb{Q}|\Pr))\, .
$$
(The measure $\mathbb{Q}$ is supposed to be Gibbsian or only ergodic if we
invoke the Ornstein-Weiss
theorem alluded to above.)
This gives an interpretation of relative entropy in terms of repetition
and waiting times.

We now turn to the analysis of fluctuations of occurrence and repetitions
of patterns.
In the sequel, $U$ is the interaction defining the
Gibbs measure $\Pr$ (see Section \ref{Gibbs} below).
The following two theorems are proved in Subsection \ref{cltproof}.

\bt\label{corollary-clt-waiting}
Let $U$ be a finite range, translation-invariant interaction, and for
$\beta$ small enough let $\Pr_\beta$ be the unique
Gibbs measure with interaction $\beta U$. There exists $\beta_0 >0$ such
that for all $\beta<\beta_0$
there exists $\theta=\theta_\beta >0$ such that
\be
\frac{\log\w_n - \E\log\w_n}{ n^{\frac{d}{2}}}
\to \mathcal{N}(0,\theta^2)\ ,\ \text{as}\ n\to\infty,\ \text{in }
\Pr_\beta\times \Pr_\beta \ \text{ distribution}\ .
\ee
where $\mathcal{N}(0,\theta^2)$ denotes the normal law with
mean zero and variance $\theta^2$, which is equal to
\be
\frac{d^2}{dq^2} \left(P((1-q)\beta U)\right)\big|_{q=0}\, .
\ee
\et

\bt\label{corollary-clt}
Let $U$ be a finite range, translation-invariant interaction,
and for $\beta$ small enough let $\Pr_\beta$ be the unique
Gibbs measure with interaction $\beta U$. There exists $\beta_0 >0$ such
that for all $\beta<\beta_0$
there exists $\theta=\theta_\beta >0$ (the same as in the previous
theorem) such that
\be
\frac{ \log\r_n- \E\log \r_n}{ n^{\frac{d}{2}}}
\to \mathcal{N}(0,\theta^2)\ ,\ \textup{a.s.}\ n\to\infty,\ \text{in }
\Pr_\beta \ \text{distribution}\, .
\ee
\et

\br
From the proof of the previous theorem it follows that
one can replace the measure $\Pr_\beta\times \Pr_\beta$
by the measure ${\mathbb Q}\times \Pr_\beta$, where ${\mathbb Q}$ is any
ergodic random field,
and $s(\Pr_\beta)$ by $s({\mathbb Q})+ s({\mathbb Q}|\Pr_\beta)$.
\er

\br
The $\beta_0$ of Theorems \ref{corollary-clt} and
\ref{corollary-clt-waiting} determines the
analyticity regime of the pressure.
This is related to the regime where the high-temperature expansion is
convergent. The restriction
to finite range interactions is here for convenience only, and can be
replaced by the requirement that
the norm
\[
\| U\|= \sum_{A\ni 0} \| U(A,\cdot)\|\ \exp\left(\alpha (\textup{diam}
(A))\right) 
\]
is finite for some $\alpha >0$, see \cite{simon}.
\er

We end our corollaries with large deviation estimates. In the context of
Gibbs measures,
it is well-known that the sequence $\{-\frac{1}{n^d} \log \Pr({\mathcal
C}(\si_{C_n})):n\in\N\}$
satisfies a large deviation principle see e.g., \cite{Comets},
\cite{Olla}.
Here we shall apply the more specific large deviation result of \cite{PS}
that was already used in \cite{cgs} to establish large deviations
for $\log \r_n$ (in dimension one).

The following theorem is proved in subsection \ref{LD-proof}.
\bt\label{LD-waiting}
Let $\Pr$ be a translation-invariant Gibbs random field satisfying
\eqref{exponmix}. Then
for all $q\in\R$ the limit
\be
\mathcal{W} (q)= \lim_{n\to\infty} \frac{1}{n^d} \log \int \w_n^{qd}\
d\Pr\!\times\!\Pr
\ee
exists. Moreover, 
\be\label{w-cases}
\mathcal{W} (q)= \begin{cases} 
P\left((1-q)U\right) + (q-1)P(U), & \mbox{ for }q\ge -1,\\
P(2U) - 2P(U), &\mbox{ for }q<-1,
\end{cases}
\ee
where $P$ is the pressure defined in (\ref{pres}) below.
\et

The following theorem gives the precise consequence of Theorem
\ref{LD-waiting}
for the large deviations of $\log \w_n$ provided 
$P((1-q)U)$ is $C^1$  for all $q\ge -1$.
For this we can apply the result  of \cite{PS}.
The pressure function is $C^1$ for example in the Ising model.
In the case $P((1-q)U)$ is not differentiable everywhere on $[-1,\infty)$,
the result of \cite{PS} will give us Large Deviations for 
$u$ in  some bounded interval. 

\bt\label{psld-waiting}
Suppose $U$ is a finite range translation-invariant interaction. Then there
exists $\beta_1 >0$ such that
for $\beta\leq \beta_1$ there exists a unique Gibbs measure $\Pr_\beta$
with interaction $\beta U$, and   for all $u\ge 0$ we have
\be\label{plat-waiting}
\limn \frac{1}{n^d} \log\ (\Pr_\beta\times \Pr_\beta)\left( \frac{\log
\w_n^d}{n^d}\geq s(\Pr_\beta) +u\right)
=
\inf_{q>-1} \left\{-(s(\Pr_\beta)+u) q + \mathcal{W} (q)\right\}
\ee
and for all $u\in (0, u_0)$, $u_0= |\lim_{q\downarrow -1} 
\mathcal W'(q)-s(\Pr)|$,
\be\label{plat-waiting-bis}
\limn \frac{1}{n^d} \log\ (\Pr_\beta\times \Pr_\beta)\left( \frac{\log
\w_n^d}{n^d}\leq s(\Pr_\beta) -u\right)
=
\inf_{q>-1} \left\{-(s(\Pr_\beta)-u) q + \mathcal{W} (q)\right\}
\ee
\et 

\br
A more general version of Theorem \ref{LD-waiting} can be easily deduced
by following
the same lines as its proof: The measure $\Pr \times \Pr$ can be replaced
by the measure ${\mathbb Q}\times \Pr$ where ${\mathbb Q}$ is any Gibbsian
random field
(without any mixing assumption). Of course formula \ref{w-cases} has to be
modified: Now
${\mathcal W}(q) = P(V-q U) - P(V)+q P(V)$ for $q\geq -1$, where $V$ is
the interaction
of the Gibbs measure ${\mathbb Q}$. Accordingly, a version of Theorem
\ref{psld-waiting}
can be obtained under a differentiability  condition on $\mathcal W$.
\er

\br
Under the assumption of Theorem \ref{LD-waiting}, the sequence
$\left\{\frac{d}{n^d}\log\w_n\right\}$
satisfies a Large Deviation Principle in the sense of \cite{DZ} (Theorem
4.5.20 p. 157). 
\er

The following theorem derives from
Theorem \ref{corollary-repetition}.
Since its derivation follows verbatim  
along the  lines of \cite{cgs}, we omit the proof.

\bt\label{psld}
Suppose $U$ is a finite range, translation-invariant interaction. There
exists $\beta_1 >0$ be such that
for $\beta\leq \beta_1$ there exists a unique Gibbs measure $\Pr_\beta$
with interaction $\beta U$ and there exists $\tilde{u}>0$ such that for
all $u\in [0,\tilde{u})$
we have
\be\label{platreturn}
\limn -\frac{1}{n^d} \log \Pr_\beta\left( \frac{\log \r_n^d}{n^d}\geq
s(\Pr_\beta) +u\right)
=\mathcal I( s(\Pr_\beta)+u),
\ee
and
\be\label{plat}
\limn -\frac{1}{n^d} \log \Pr_\beta\left( \frac{\log \r_n^d}{n^d}\leq
s(\Pr_\beta) -u\right)
=\mathcal I( s(\Pr_\beta)-u),
\ee
where 
$$
\mathcal{I}( u)=\sup_{q\in\R} \left(u q - P((1-q)U)-(q-1)P(U)\right)
$$
\et
\br 
It follows from the proof of theorem \ref{psld-waiting}
that
we have the analogue
theorem for repetition times, if
we condition the measure $\Pr_\beta$ 
on good patterns, that is, patterns which are not ``badly
self-repeating'', see Definition \ref{bsr}
below. 
\er
\br
$\beta_1$ in Theorems \ref{psld-waiting} and \ref{psld} does not
necessarily coincide with the critical
inverse temperature $\beta_c$ (below which there is a unique
Gibbs measure), and is in general strictly larger than $\beta_0$ of
Theorem
\ref{corollary-clt} and \ref{corollary-clt-waiting}, see \cite{enter}.
\er

\section{Gibbsian random fields and Dobrushin uniqueness }\label{Gibbs}

For the sake of convenience
the present and next subsections are devoted to the notion of Gibbsian
random
fields and their mixing properties. More details
on this subject can be found in \cite{Geo}, \cite{Guyon}.
\bd
A translation-invariant interaction is a function
\be\label{U}
U:\s\times\Omega\rightarrow\R,
\ee
such that the following conditions are satisfied:
\begin{enumerate}
\item
$U(A,\sigma)$ depends on $\sigma ({\bf x})$,
with ${\bf x}\in A$ only.
\item
{\it Translation invariance:} 
\be\label{UAx}
U(A+{\bf x},\tau_{-{\bf x}}\sigma)= U(A,\sigma ) \qquad \forall A\in\s,
{\bf x}\in\Z^d, \sigma\in \Omega.
\ee
\item
{\it Uniform summability:}
\be\label{sumA}
\sum_{A\ni 0} \sup_{\sigma\in \Omega} | U(A,\sigma )| <\infty \ .
\ee
\end{enumerate}
\ed

An interaction $U$ is called {\it finite-range} if there exists an $R>0$
such that
$U(A,\sigma )=0$ for all $A\in \s$ with $\mbox{diam}(A)>R$. 

The set of all such interactions is denoted by ${\mathcal U}$. Mostly
we will give examples of Gibbs measures satisfying our mixing
conditions with interactions $U\in {\mathcal U}$. This can be generalized
easily to interactions such that
\[
\| U\|_\alpha= \sum_{A\ni 0} \| U(A,\cdot)\|\ \exp\left(\alpha (\textup{diam}
(A))\right) 
\]
is finite for some $\alpha >0$.

For $U\in {\mathcal U}$, $\zeta\in\Omega$, $\Lambda\in\s$, we define the
finite-volume Hamiltonian
with boundary condition $\zeta$ as
\be\label{Hzeta}
H^\zeta_\Lambda (\sigma )= \sum_{A\cap\Lambda\not=\emptyset}
U(A,\sigma_\Lambda\zeta_{\Lambda^c})\, .
\ee
Corresponding to the Hamiltonian in (\ref{Hzeta}) we have the
finite-volume Gibbs measures $\Pr^{U,\zeta}_\Lambda$,
$\Lambda\in\s$, defined on $\Omega$ by
\be\label{finvol}
\int f(\xi)\ d\Pr_\Lambda^{U,\zeta} (\xi)
= \sum_{\sigma_\Lambda\in\Omega_\Lambda}
f(\sigma_\Lambda\zeta_{\Lambda^c})\
e^{-H^\zeta_\Lambda (\sigma)}/ Z^\zeta_\Lambda\, ,
\ee
where $f$ is any continuous function and $Z_\Lambda^\zeta$ denotes the
partition function normalizing $\Pr^{U,\zeta}_\Lambda$ to
a probability measure. Because of the uniform summability condition,
(\ref{sumA}) the objects $H_\Lambda^\zeta$ and $\Pr^{U,\zeta}_\Lambda$
are continuous as functions of the boundary condition $\zeta$.

For a probability measure $\Pr$ on $\Omega$, we denote by
$\Pr^\zeta_\Lambda$ the
conditional probability distribution of $\sigma ({\bf x}), {\bf x}\in
\Lambda$, given
$\sigma_{\Lambda^c}=\zeta_{\Lambda^c}$. Of course, this object is only
defined on
a set of $\Pr$-measure one. For $\Lambda\in\s, \Gamma\in \s$ and
$\Lambda\subset\Gamma$,
we denote by $\Pr_\Gamma(\sigma_\Lambda|\zeta)$ the conditional
probability to find
$\sigma_\Lambda$ inside $\Lambda$, given that $\zeta$ occurs in
$\Gamma\setminus\Lambda$.

For $U\in {\mathcal U}$, we call $\Pr$ a Gibbs measure with interaction
$U$ if its conditional
probabilities coincide with the ones prescribed by (\ref{finvol}), i.e.,
if
\be\label{gibbsdef}
\Pr^\zeta_\Lambda = \Pr^{U,\zeta}_\Lambda \qquad \Pr-a.s. \qquad
\Lambda\in\s,
\zeta\in\Omega.
\ee
We denote by $\mathcal{G} (U)$ the set of all translation invariant
Gibbs measures with interaction $U$. For
any $U\in {\mathcal U}$, $\gee (U)$ is a non-empty compact convex set.
In this paper we will in fact restrict ourselves to interactions with a
unique Gibbs measure.

A basic example is the ferromagnetic Ising model, where
$U(\{{\bf x},{\bf y}\},\si) = -\beta J\si({\bf x})\si({\bf y})$ if $|{\bf
x}-{\bf y}|=1$, 
$U(\{{\bf x}\},\si) =-h\beta\si ({\bf x})$. Here $\beta\in (0,\infty)$
represents the inverse temperature,
$J>0$ the coupling strength, and $h$ the external magnetic field.

\bigskip

We turn to the mixing properties of Gibbs random fields.
For an interaction $U\in{\mathcal U}$, the Dobrushin matrix is given by
$$
\gamma_{{\bf xy}} (U) =
\frac{1}{2}\sup\left\{
|\ \Pr^{U,\zeta}_{\{{\bf x}\}} (\alpha ) -\Pr^{U,\xi}_{\{{\bf x}\}}
(\alpha)|:\
\zeta,\xi\in\Omega, \zeta_{\Z^d\setminus\{{\bf y}\}}
=\xi_{\Z^d\setminus \{{\bf y}\}}\ , \, \alpha\in \mathcal{A}\right\}\ .
$$
The matrix $\gamma$ measures the dependence
of changing the spin at site ${\bf y}$ on the conditional
probability at site ${\bf x}$.

\bd\label{dobrushin}
The interaction $U$ is said to satisfy the Dobrushin
uniqueness condition
if
\be\label{dobuniq}
\sup_{{\bf x}\in\Z^d}\sum_{{\bf y}\in\Z^d} \gamma_{{\bf xy}} (U) <1\, .
\ee
\ed
The following result is proved in \cite{Geo}, see
also \cite{Guyon}, theorem 2.1.3, p. 52.
\bt
Let $U\in{\mathcal U}$ be a finite range interaction.
Under the condition (\ref{dobuniq}), there is a
unique Gibbs measure $\Pr\in\gee (U)$, and
this $\Pr$ is non-uniformly exponentially $\varphi$-mixing,
i.e., it satisfies the mixing property (\ref{exponmix}).
\et

Examples for which (\ref{dobuniq}) is satisfied are:
\begin{enumerate}
\item
The so-called high-temperature region where.
$U\in{\mathcal U}$ is such that
\be\label{DobU}
\sup_{{\bf x}\in\Z^d} \sum_{A\ni {\bf x}}
(|A|-1)\sup_{\sigma, \sigma'\in\Omega}|U(A,\sigma)-U(A,\sigma')| < 2.
\ee
Inequality (\ref{DobU}) implies the Dobrushin uniqueness condition 
(\ref{dobuniq})
(see \cite{Geo}, p.\ 143, Proposition 8.8).
In particular, it implies that $|\gee (U)|=1$ (i.e., no phase transition).
Note that
it is independent of the ``single-site part'' of the interaction, i.e., of
the
interactions $U(\{{\bf x}\},\sigma)$.
For any finite range potential $U$ there exists $\beta_c$ such that $\beta
U$ 
satisfies \eqref{DobU} for all $\beta<\beta_c$.
For the Ising model in $\Z^2$, much more is known: the mixing property
\eqref{exponmix}
holds for {\em any} $\beta <\beta_c$
(see e.g. \cite{dobshlos}).
\item Low temperature regime
for an interaction with unique ground state, e.g.,
the Ising model in a homogeneous
magnetic field and sufficiently large $\beta$. See \cite{Guyon} example
(2.1.5)
\item
Interactions in
a large external field. See \cite{Guyon}, example
(2.1.4). For the Ising model in two dimensions this means
that the field $h$ should satisfy
\[
|h| > 4\beta + \log (8 \beta)\, .
\]
\end{enumerate}
\br
The Dobrushin uniqueness condition is
not a necessary condition for
the mixing property
\ref{mixdef}. More general versions, known
as ``Dobsruhin-Shlosman" conditions exist,
see e.g., \cite{martinelli} for more details on
general finite size conditions ensuring NUEM.
\er

We now recall some basic facts on entropy and relative entropy (or
Kullback-Leibler information).
We use the following shorthand to ease notation~:
$$
\sum_{\C_n}=\sum_{\C(A_n):A_n\in\Omega_{C_n}}\;.
$$

The entropy $s(\Pr)$ of $\Pr$ is defined as
$$
s(\Pr)=\lim_{n\to\infty}-\frac1{n^d}\sum_{\C_n}\Pr(\C_n)\log\Pr(\C_n)\, .
$$

The relative entropy $s(\mathbb{Q}|\Pr)$ of a stationary random field
$\mathbb{Q}$ with respect to a Gibbsian
random field $\Pr$ is
$$
s(\mathbb{Q}|\Pr)=\lim_{n\to\infty}\frac1{n^d}\sum_{\C_n}\mathbb{Q}(\C_n)\log\frac{\mathbb{Q}(\C_n)}{\Pr(\C_n)}
$$
In terms of the interaction $U$ of $\Pr$ the relative entropy is
\[
s(\mathbb{Q}|\Pr)=P(U)+\int f_U\ d\mathbb{Q}-s(\mathbb{Q}),
\] 
where
\[
f_U (\si) =\sum_{A\ni 0} \frac{U(A,\si)}{|A|}
\]
and $P(U)$ is the pressure of $U$, which defined as follows
\be\label{pres}
P(U) =\lim_{n \to\infty}\frac{1}{n^d}\log Z_{C_n}\, ,
\ee
where
\[
Z_{C_n} = \sum_{\si_{C_n}\in\Omega_{C_n}} \exp(-\sum_{A\subset C_n}
U(A,\si))
\]
is the partition function with the free boundary conditions.


\bp\label{asre}
Let $\Pr$ be a Gibbs random field and $\mathbb{Q}$ be an ergodic random
field. Then
\[
\lim_{n\to\infty}\frac1{n^d}\log \frac{\mathbb{Q}({\mathcal
C}(\si_{C_n}))}{\Pr({\mathcal C}(\si_{C_n}))}
=s(\mathbb{Q}|\Pr)
\]
for $\mathbb{Q}$-almost every $\si$.
\ep

\bpr
The proof is simple, but since we did not find it in the literature,
we give it here for the sake of completeness. Write
\[
g_n (\si) \sim h_n (\si)
\]
if
\[
\lim_{n\to\infty}\frac{1}{n^d}\sup_\si |g_n(\si)-h_n (\si)|=0\, .
\]
Let $U$ be the potential of the Gibbsian field $\Pr$.
Then we have
\[
\log \Pr({\mathcal C}(\si_{C_n})) \sim -\sum_{i\in C_n} \tau_i f_U
(\si) -\log Z_{C_n}\, .
\]

Therefore, by ergodicity of $\mathbb{Q}$
\[
\frac{1}{n^d} \log \Pr({\mathcal C}(\si_{C_n}))
\]
converges $\mathbb{Q}$-a.s. to
\[
-\int f_U d\mathbb{Q} - P(U)\, .
\]
By the Shannon-Mc Millan-Breiman theorem \cite{follmer,thouvenot}
\[
\frac{1}{n^d}\log {\mathbb Q}({\mathcal C}(\si_{C_n}))
\]
converges $\mathbb{Q}$-a.s. to $-s({\mathbb Q})$.
Hence
the difference 
\[
\frac{1}{n^d}(\log \mathbb{Q}({\mathcal C}(\si_{C_n})) -\log \Pr({\mathcal
C}(\si_{C_n})))
\]
converges $\mathbb{Q}$-a.s. to
\[
P(U)-\bigl(s(\mathbb{Q})-\int f_U\ d\mathbb{Q} \bigr)
\]
which is equal to $s(\mathbb{Q}|\Pr)$ by the Gibbs variational principle, see\cite{Geo}. \epr

A standard property of Gibbs measures which we will use
often is the following: there exist positive constants $C,c,C',c'$ such that
\be\label{gibbsineq}
 C'\ e^{ -c' n^d}\le \Pr( \mathcal C_n ) \le C\ e^{-cn^d}
\ee
for every cylinder $\mathcal C_n$ supported on $C_n$.

            \section{Proof of Theorem
\protect{\ref{thm1}}}\label{mainproof}

To ease notation, we will write $\Pr(A)$ instead of $\Pr({\mathcal C}(A))$
where $A=A_n$
is an $n$-pattern.

\subsection{Preliminary results}

In this section we prove Theorem \ref{thm1}.
We follow the approach of \cite{abadi}.

For $V\in\s$, $\si\in\Omega$ and $A=A_n$ an $n$-pattern we say that
``$A$ is present in $V$'', and write $A\prec V$, for the configuration
$\si$
if there exists ${\bf x}=(x_1, x_2, ...,x_d)\in \Z^d$ such that $W:= {\bf
x}+C_n\subset V$
and $(\tau_{{\bf x}}\si)_W= A$. By abusing notation, we will
write $\Pr(A\prec V)$ for the probability of that event.

\bl\label{volume_estim} 
Let $V$ be a finite subset of $\Zd$, and let $A=A_n$ be a 
$n$-pattern. Then 
$$
\Pr(A\prec V) \le |V|\ \Pr(A).
$$
\el
\begin{proof}$\displaystyle
\Pr(A\prec V)\le \sum_{\mathbf x\in V} \Pr\bigl(\bigl\{
\sigma:\, \sigma|_{\mathbf x+C_n}=A\bigr\}\bigr)
=\sum_{\mathbf x\in V} \Pr(A) = |V|\ \Pr(A).
$
\end{proof}

For every $k\in\N$ define
$$
 N_k^A(\si) = \sum_{\begin{subarray}{c} 
{\mathbf x}\in\Z^d\\
0\le x_i\le k
\end{subarray}}
\1\{ \tau_{{\bf x}}(\si)_{C_n}=A\}.
$$
Then the following events coincide:
\begin{equation}\label{equalevents}
\{ \t_A\le k\} = \{ N_k^A \ge 1\}.
\end{equation}
Moreover,
$$
  \E N_k^A =  (k+1)^d \Pr(A).
$$

\bl[Second moment estimate]\label{secondmoment} 
Consider a non-uniformly exponentially $\phi$-mixing Gibbsian random
field. Then there exists $\delta>0$ such 
that for every $n, k\in \N$, and every $\Delta>2n$ one has
$$
\E (N_k^A)^2 \le
(k+1)^d\Pr(A) \left(1 + e^{-\delta n}\Delta^d +
(k+1)^{d}\Pr(A)
+ (k+1)^{d} n^d\phi(\Delta-2n)\right).
$$
\el
\begin{proof}
Define $C(\mathbf x,n)= \mathbf x+C_n$.
We have to estimate the following expression
\begin{equation}
\E (N_k^A)^2  = \sum_{\begin{subarray}{c} 
{\mathbf x}\in\Z^d\\
0\le x_i\le k
\end{subarray}}  
\sum_{\begin{subarray}{c} 
{\mathbf y}\in\Z^d\\
0\le y_i\le k 
\end{subarray}}\Pr( \si_{C({\bf x},n)}= \si_{C({\bf y},n)}=A).
\end{equation}
We split the above double sum into the three following sums
\begin{equation}
I_1 = \sum_{\mathbf x=\mathbf y},\quad 
I_2 =\sum_{
\begin{subarray}{c}
\mathbf x\ne \mathbf y\\ \nonumber
|\mathbf x-\mathbf y|\le \Delta
\end{subarray}},
\quad
I_3  =\sum_{
\begin{subarray}{c}
\mathbf x\ne \mathbf y\\ \nonumber
|\mathbf x-\mathbf y|> \Delta
\end{subarray}}
\end{equation}
Let us proceed with each of the sums separately. For $I_1$ one obviously
has
$$
 I_1 = \sum_{\begin{subarray}{c} 
{\mathbf x}\in\Z^d\\
0\le x_i\le k
\end{subarray}}  \Pr(A) = (k+1)^d \Pr(A).
$$
To estimate $I_2$ we use that for any Gibbsian random field
there exists a constant $\delta>0$ such
that for any finite volume V, and any configuration $\sigma$ and $\eta$,
the conditional probability of observing $\sigma$ on $V$, given $\eta$
outside of $V$,
can be estimated as follows \cite{Geo}  
$$\Pr( \sigma_V\, | \eta_{V^c}) \le \exp( -\delta |V|).
$$
Therefore
$$
\aligned
I_2 &= \sum_{\begin{subarray}{c}
\mathbf x\ne\mathbf y\\
|\mathbf x- \mathbf y|\le \Delta
\end{subarray}} \Pr( \si_{C({\bf x},n)} =\si_{C({\bf y},n)} =A) 
=\sum_{\begin{subarray}{c}
\mathbf x\ne\mathbf y\\
|\mathbf x- \mathbf y|\le \Delta
\end{subarray}} 
\Pr\bigl( \si_{C({\bf x},n)}=A \,\bigl|\, \si_{C({\bf y},n)}
 =A\bigr.\bigr)\ \Pr(A)\\
&\le  \sum_{\begin{subarray}{c}
\mathbf x\ne\mathbf y\\
|\mathbf x- \mathbf y|\le \Delta
\end{subarray}}
\Pr(A) \exp\bigl( -\delta\ | C({\bf x},n)\setminus C({\bf y},n)|)\ .
\endaligned
$$
To complete the estimate, it is sufficient to observe that since $\mathbf
x\ne\mathbf y$,
the volume of the set  $C({\bf x},n)\setminus C({\bf y},n)$ is at least n.
Hence
$$
  I_2 \le (k+1)^d\Delta^d \exp(-\delta n)\ \Pr(A).
$$
Finally, using the mixing condition \eqref{exponmix}, for $I_3$ we obtain
$$
\aligned 
I_3 &= \sum_{\begin{subarray}{c}
\mathbf x\ne\mathbf y\\
|\mathbf x- \mathbf y|>\Delta
\end{subarray}} \Pr( \si_{C({\bf x},n)} =\si_{C({\bf y},n)} =A)
\le \sum_{\begin{subarray}{c}
\mathbf x\ne\mathbf y\\
|\mathbf x- \mathbf y|>\Delta
\end{subarray}}\Bigl( \Pr(A)+ n^d\phi(\Delta-2n)\Bigr)\ \Pr(A)\\
&\le (k+1)^{2d}\ \Pr(A) \bigl( 
\Pr(A) +  n^d\phi(\Delta-2n)
\bigr).
\endaligned
$$ 
Combining all the estimates together we obtain the statement of the lemma. 
\end{proof}

\bl[The parameter]\label{pospar}
There exist strictly positive constants $\Lambda_1,
\Lambda_2$ such that
for any integer  $t$  with $t\Pr(A)\le 1/2$, one has 
$$
    \Lambda_1\le \lambda_{A,t} := -\frac { \log\Pr( \t_A>t^{1/d})}{t
\Pr(A)}
\le \Lambda_2.
$$
\el
\begin{proof} Taking into account (\ref{equalevents}) and the 
Cauchy-Schwartz inequality we obtain
\be\label{bruu} 
\Pr(\t_A\le k) \ge \frac { (\E N_k^A)^2}{ \E (N_k^A)^2}.
\ee
We apply the basic inequalities
\begin{equation}\label{elementary}
  \frac \kappa 2 \le 1-e^{-\kappa} \le \kappa\ ,
\end{equation}
where the left inequality is valid for all $\kappa \in [0,1]$,
and the right inequality is true for $\kappa\ge 0$.
Let now $\kappa = -\log\Pr( \t_A>t^{1/d})$. Then, using 
lemma \ref{secondmoment} and (\ref{bruu}), we conclude 
$$
\aligned
\frac {-\log\Pr( \t_A>t^{1/d})}{ t\Pr(A)} & \ge 
\frac {\Pr( \t_A\le t^{1/d})} { t\Pr(A)} \\
&\ge \frac {1}{1+e^{-\delta n}\Delta^d+ (t+1)\Pr(A)+ (t+1) n^d
\phi(\Delta-2n)}\\
&\ge \frac {1}{1+c_1+3/2 + c_2}=:\Lambda_1,
\endaligned
$$
where we have chosen $\Delta = n^{d+1}$, 
$$ 
c_1 = \sum_{n\in \N} e^{-\delta n} n^{d(d+1)} <\infty, \ \text{and}\ 
c_2 = \sup_{n\in \N}\left\{ (t+1)\, n^d \phi(\Delta - 2n)\right\} 
$$
We have to show that $c_2$ is finite. Indeed, since for  a Gibbs random
field $\Pr$ there exist $c',C'>0$ such that
$$
\Pr(A) \ge C'\exp( -c'n^d)
$$
for every $n$-pattern $A$; $t$ has been chosen such that $t\Pr(A)<1/2$, we have
$$
c_2\le  \sup_{n\in\N} \left\{ \left[\frac 1{2C'}\exp( c'n^d)+1\right]n^dC_1\exp\bigl(-C_2(n^{d+1}-2n)\bigr)\right\} <\infty,
$$
where we have used the mixing condition (\ref{exponmix}).

For the upper bound, we use 
(\ref{elementary}) again, but first we have to check that 
$$\kappa= - {\log\Pr( \t_A>t^{1/d})}\in [0,1].$$
Indeed, since $t<(2\Pr(A))^{-1}$, by Lemma \ref{volume_estim}  we have
$$
\Pr( \t_A>t^{1/d})\ge \Pr \left(\t_A > \frac {1}{(2\Pr(A))^{1/d}}\right)
=1 - \Pr \left(\t_A \le \frac {1}{(2\Pr(A))^{1/d}}\right)
\ge 1 - \frac {\Pr(A)}{2\Pr(A)}=\frac 12.
$$ 
Hence, $\kappa\le \log(2) < 1$, therefore 
$\kappa \le 2(1-e^{-\kappa})$, which means  
$$
  - {\log\Pr( \t_A>t^{1/d})}\le 
2{\Pr( \t_A\le t^{1/d})}\le 2 t\ \Pr(A)\leq 1
$$
where we have used Lemma \ref{volume_estim} for the second inequality.
Hence, we can choose $\Lambda_2 = 2$. This finishes the proof.\end{proof}

For positive numbers $x_{A_n}, y_{A_n}$ depending on the $n$-pattern
$A_n$ we write $x_{A_n} \sim y_{A_n}$ if
$$
\lim_{n\to\infty} \frac{x_{A_n}}{y_{A_n}} = 1.
$$
For a positive integer $t_A$  we set
$C(t_A)= [0,t_A]^d \cap {\mathbb Z}^d$.
For a subset $V\subset \Z^d$ let 
$A\nprec V$ be the event that the $n$-pattern $A$ cannot
be found in $V$. (See above for the definition of $A\prec V$.)

The following lemma is crucial and gives the factorization
property of the exponential distribution, i.e.,
the fact that asymptotically
\[
\Pr (t_{A_n} > (t+s)/\Pr (A_n))\sim 
\Pr (t_{A_n} > t/\Pr (A_n))\Pr (t_{A_n} > s/\Pr (A_n))
\]
where the accuracy of the approximation
marked $\sim$ is spelled out in detail. The idea is
that the event of non-occurrence of the pattern
in a cube of size $O(1/\Pr (A_n))$ can be viewed
as the non-occurrence of the pattern in
many sub-cubes of volume $k_n$, where $n^d<< k_n<< 1/\Pr (A_n)$.
These sub-cubes will be separated by corridors
of width $\Delta$, where $\Delta$ is such that
the pattern occurs with very small probability
in the corridor, and on the other hand
the mixing can be used to decouple the events
of non-occurrence in different sub-cubes.

Our choice for the volume of the sub-cubes
will be $k_n = O(\Pr(A_n)^{-\theta})$, with $\theta\in (0,1)$ and the corridors
will have width $\Delta_n = O(n^k)$ with $k$ big enough for
the mixing to work well.

This explains the choices in the statement of the following lemma. 
\begin{lemma}[Iteration Lemma]\label{iterationlemma}
Let $A=A_n$ be a $n$-pattern and $t_A$ be such that
$t_A^d=[\Pr(A)^{-\vartheta}]$, where $[\cdot]$
denotes the integer part, and
$\vartheta\in (0,1)$. For $i=1,\ldots k$,
let
$C_i (t_A)$ denote any collection of $k$ disjoints
cubes of the form ${\mathbf x}_i+ C(t_A)$.
Then, for $n$ large enough, 
there exists $\delta\in(0,1)$, which depends only on the measure $\Pr$, 
such that the following inequality holds for all $k$:

\be\label{ineq}
\biggl|\Pr\left( A\nprec\bigcup_{i=1}^k C_i (t_A) \right)
-
\Pr \left( A\nprec C(t_A)\right)^k \biggr|
\leq 
k\Pr(A)^{\eta} 
\Bigl( 
\Pr\left( A\nprec C(t_A)\right)+\Pr(A)^{\eta}
\Bigr)^{k},
\ee
where $\eta=\left( 1- \vartheta (d-1)/d \right) (1-\delta).$
\end{lemma}
\bpr
We will prove that
\be\label{ineqtoprove}
\Pr\left(  A\nprec\bigcup_{i=1}^k C_i (t_A)\right)
\leq 
\Pr\left( A\nprec C(t_A)\right)^k 
+ 
\Pr(A)^{\eta} \ k \
\left( 
\Pr\left( A\nprec C(t_A)\right)+\Pr(A)^{\eta}
\right)^{k}
\ee
The inequality
\be \nonumber
\Pr\left(A\nprec\bigcup_{i=1}^k C_i (t_A)\right)
\geq \Pr\left(A\nprec C(t_A)\right)^k 
   - 
\Pr(A)^{\eta} \ k \
\left( 
\Pr\left( A\nprec C(t_A)\right)+\Pr(A)^{\eta}
\right)^{k}
\ee
is derived analogously.

For any positive integer $z$, we write ${\bf z}=(z,\dots,z)\!\in\!\Z^d$.
We denote by $C_i^z(t_A)$ the cube ${\mathbf x_i}+{\bf z}+C(t_A-2 z)$
and by $C^z(t_A)$ the cube $C(t_A -2z)$.
For any positive integer $\Delta < 2t_A$,
we consider the difference
\beq
&& 
\left\vert
\Pr\left( A\nprec\bigcup_{i=1}^k C_i (t_A) \right) -
\Pr\left( A\nprec C_1^\Delta(t_A)\cup\bigcup_{i=2}^k C_i (t_A)  \right)
\right\vert\nonumber\\
&=&
\Pr\left( \big(A\nprec C_1^\Delta(t_A)\cup\bigcup_{i=2}^k C_i (t_A)\big) 
            \ \cap \ \big(A\prec C_1(t_A)\backslash C_1^\Delta(t_A)\big)
\right) \nonumber\\
&\le&
\Pr\left( \big(A\nprec \bigcup_{i=1}^k C_i^{2\Delta} (t_A)\big) 
            \ \cap \ \big( A\prec C_1(t_A)\backslash C_1^\Delta(t_A)\big)
\right). \nonumber
\eeq
Iterating the mixing property \eqref{mixfun} and using Lemma
\ref{volume_estim}, we bound the last term by
\beq \nonumber
2 d \Delta t_A^{d-1}\Pr(A) \
\left(
\Pr\left( A\nprec C^{2\Delta} (t_A)  \right) + 
          \phi(\Delta)\ |C(t_A)| \
\right)^k \  .
\eeq

On the other hand
$$
\aligned
\Bigl\vert
\Pr\Bigl( A\nprec C_1^\Delta(t_A) & \cup\bigcup_{i=2}^k C_i (t_A)  \Bigr)-
\Pr\left( A\nprec C_1^\Delta(t_A)  \right) 
\Pr\Bigl( A\nprec \bigcup_{i=2}^k C_i(t_A)  \Bigr)
\Bigr\vert\\
&\leq
\phi(\Delta)
\left\vert C_1^\Delta(t_A) \right\vert 
\ \Pr\left( A\nprec \bigcup_{i=2}^k C_i(t_A)  \right)\\
&\leq
\phi(\Delta) t_A^{d}
\left(
\Pr\left( A\nprec C^{\Delta} (t_A)  \right) + \phi(\Delta)\ |C(t_A)| \
\right)^{k-1}.
\endaligned
$$

Put
$$\aligned
\epsi_1 &= \epsi_1(A,t_A,\Delta)= \varphi (\Delta) \ t_A^d,\\
\epsi_2 &= \epsi_2 (A,t_A,\Delta) =2 d \Delta t_A^{d-1} \Pr(A),\\
\epsi_{\phantom{1}}   &= C (\epsi_1 + \epsi_2),
\endaligned
$$
where $C$ is a positive constant to be defined later on. Put also
\be \nonumber
\alpha_{k-j} = \Pr\left(A\nprec \bigcup_{i=j+1}^{k} C_i(t_A)\right),
\ee
\be \nonumber
\alpha_{k-j}^z= \Pr\left(A\nprec \bigcup_{i=j+1}^{k} C_i^z(t_A)\right)
.
\ee
We obtain the recursion
\be \nonumber
\alpha_k 
\leq (\epsi_1 + \epsi_2) (\alpha_1^{2\Delta} +\epsi_1)^{k-1}
+\alpha_1^\Delta\alpha_{k-1}\ , 
\ee
which upon iteration leads to
\be \nonumber
\alpha_k \leq (\epsi_1 + \epsi_2) \ k \ (\alpha_1^{2\Delta}+\epsi_1)^{k-1}
+ (\alpha_1^\Delta)^k .
\ee
We choose $C$ such that, $\alpha_1^{2\Delta}+\epsi_1 \le \alpha_1+\epsi$, 
so we have
\be \nonumber
\alpha_k \le \epsi \ k \ (\alpha_1+\epsi)^{k-1} + (\alpha_1+\epsi)^k .
\ee
Now we use the following simple inequality:
for $0< x\leq y <1$ and $N$ any positive integer
$$ 
y^N-x^N = (y-x) (x^{N-1} +x^{N-2}y + \ldots + y^{N-1})
\leq 
(y-x) N y^{N-1}
$$
to obtain
\be\label{lastineqrequired}
\alpha_k -\alpha_1^k 
\leq 2 \ \epsi \ k \ (\alpha_1^\Delta +\epsi_1)^{k-1} 
\leq 4 \ \epsi \ k \ (\alpha_1 + \epsi)^{k}
\ .
\ee
Choose $\Delta = n^{d+1}$.
Since $t_A^d\sim \Pr(A)^{-\vartheta}$ for
some $\vartheta\in (0,1)$, and since we have exponential $\phi$-mixing
\eqref{exponmix}
we obtain for the ``error terms'' $\epsi_1,\epsi_2$:
\beq
\epsi_1 &\sim & e^{-c_1 n^{d+1}} e^{c_2 \vartheta n^d} ,
\nonumber\\
\epsi_2 &\sim & 2 d n^{d+1} \Pr(A)^{-\vartheta\frac{d-1}{d}} \Pr(A)
\nonumber .
\eeq
This yields
\be
\epsi\leq
 \Pr(A)^{\left( 1- \vartheta \frac{d-1}{d} \right) (1-\delta) }, \nonumber
\ee
which together with (\ref{lastineqrequired}) implies (\ref{ineqtoprove}).
\epr

\subsection{Proof of Theorem \protect{\ref{thm1}}}

Let $t>0$, and put $t = kf_A+r$, where  
$f_A = [1/(\Pr(A))^\gamma]$ ($\gamma\in (0,1)$, $[\cdot ]$ denotes integer
part), $k$ is an integer
and $r<f_A$.
Put $t_A'= k[f_A]$,  $t_A''=(k+1)[f_A]$,
Without loss of generality we assume that the size $n$ of a $n$-pattern
$A$ is sufficiently
large, so $ f_A \Pr(A) \sim\Pr(A)^{1-\gamma} <1/2$. We remind that for 
$n$-patterns, Gibbs fields admit uniform estimates $\Pr(A)\le \exp(-cn^d)$
for some
$c>0$.
Now, recall from Lemma \ref{pospar} that
\begin{equation}\label{chooselambda}
\lambda_{A} = -\frac {\log \Pr(\t_A>(f_A)^{1/d})}{f_A\Pr(A)}
\in[\Lambda_1,\Lambda_2]
\end{equation}
for some positive constants $\Lambda_1,\Lambda_2$. 
We also define
\begin{equation} 
\tilde \lambda_{A} = 
-\frac{\log \left(\Pr(\t_A>(f_A)^{1/d})+ \Pr(A) \right)}{f_A\Pr(A) } .
\nonumber
\end{equation}
It is not difficult to see that $\tilde
\lambda_{A}\in[\Lambda_1/2,\Lambda_2]$,
for $n$ large enough.

Since $t_A' \le t \le t_A''$, one obviously has
$$
 \Pr(\t_A>t^{1/d})- \exp(-\lambda_{A}\Pr(A)\ t)\ge 
 \Pr(\t_A>(t_A'')^{1/d}) -\exp(-\lambda_{A}\Pr(A)\ t_A'),
$$
and
$$
 \Pr(\t_A>t^{1/d})-\exp(-\lambda_{A}\Pr(A) t)\le 
 \Pr(\t_A>(t_A')^{1/d}) -\exp(-\lambda_{A}\Pr(A)\ t_A'').
$$
Now, 
$$
\aligned 
|\Pr(\t_A>(t_A')^{1/d}) -\exp(-\lambda_{A}\Pr(A) t_A'')|&\le
|\Pr(\t_A>(t_A')^{1/d}\ ) -\Pr( \t_A>(f_A)^{1/d}\ )^k|\\
&+|\Pr( \t_A>(f_A)^{1/d}\ )^k-\exp( -\lambda_{A}\Pr(A)\ t_A')|\\
&+| \exp( -\lambda_{A}\Pr(A)\ t_A')-\exp( -\lambda_{A}\Pr(A)\ t_A'')|
\endaligned
$$
By Lemma \ref{iterationlemma},
\beq
|\Pr(\t_A>(t_A')^{1/d}) -\Pr( \t_A>(f_A)^{1/d} )^k|
&\le& 
\Pr(A)^{ \gamma (1-\delta)/d} \ k \
\left( \Pr( \t_A>(f_A)^{1/d} )+\Pr(A)^{1-\gamma}\right)^{k}. \nonumber\\
&=&
\Pr(A)^{ \gamma (1-\delta)/ d}\
t\ \Pr(A)  \exp( - \tilde\lambda_A \Pr(A)\ t_A')   \nonumber\\
&\le&
\Pr(A)^{ \gamma (1-\delta)/ d}\ 
t\ \Pr(A) \exp(- C_1 \Pr(A)\ t )   \nonumber .
\eeq
By the choice of $\lambda_A$ (\ref{chooselambda}) , and since $t_A' = k
f_A$,
$$
\Pr( \t_A>(f_A)^{1/d}\ )^k = \exp( - \lambda_A kf_A\Pr(A)) = \exp( -
\lambda_A t_A'\Pr(A)). 
$$
Finally,
\beq
    | \exp( -\lambda_{A}\Pr(A)\ t_A')-\exp( -\lambda_{A}\Pr(A)\
t_A'')|
&\le& \lambda_A\ \Pr(A) (t_A''-t_A') \exp( -\lambda_{A}\Pr(A)\
t_A')\nonumber\\
&\le& \Lambda_2\ \Pr(A) f_{\!A} \exp( -\Lambda_1\Pr(A)\ t_A')\nonumber\\
&\le& C_2\ \Pr(A)^{1-\gamma} \exp( -C_3 \Pr(A)\ t). \nonumber
\eeq
The lower estimate is obtained in a similar way. This finishes the proof.

\br
Notice that in the iteration lemma it is
not used that $A_n$ is a pattern. Therefore,
this lemma can be generalized to arbitrary
measurable events $E_n\in \fe_{C_{k_n}}$, where
$k_n\lll 1/\Pr(E_n)$. The second
moment estimate however uses that $A_n$ is
a pattern. Therefore Theorem 2.6  can
be generalized as follows. Let $E_n\in\fe_{C_{k_n}}$,
where $|C_{k_n}|= O(n^\alpha)$, and $\Pr(E_n ) = O(e^{-cn^d})$.
Suppose furthermore that
\be\label{secmoncon}
\limsup_{n\to\infty} \sum_{0<|x|\leq n^\alpha} 
\frac{\Pr(E_n\cap\theta_x E_n)}{\Pr (E_n)} <\infty
\ee
then (\ref{main}) holds for the
occurrence time $\t_{E_n}$.
Condition (\ref{secmoncon}) takes care of
the second moment estimate.
\er

\section{Proof of the other theorems}

\subsection{Proof of Theorem
\protect{\ref{corollary-repetition}}}\label{proof-coro-repet}
We start with a lemma on ``badly self-repeating'' patterns.
\bd
A pattern $A_n$ is called badly self-repeating  if
there exists ${\bf x}$, $0<|{\bf x}|\leq n/2$, such that
\[
\tau_{{\bf x}} {\mathcal C}(A_n)\cap{\mathcal C}(A_n)\not=\emptyset
\]
\ed
Correspondingly, a cylinder is called bad if it is of the
form ${\mathcal C} (A_n)$ with $A_n$ badly self-repeating.
The union of bad $n$-cylinders is denoted by ${\mathcal B}_n$.

\begin{lemma}[Conditioning on the initial pattern]\label{shortret}
Let $A=A_n$ be a ``good'' pattern, that is, not a badly self-repeating
pattern. 
Let $t_A$ be such that $t_A^d\sim \Pr(A)^{-\vartheta}$, where
$\vartheta\in (0,1)$.
Then there exist positive constants $b_1,b_2$ such that for 
all integers $n\geq 1$, one has
$$
\left|
\Pr\left( A\nprec C(t_A)\backslash C_n\ \big| \ A\prec C_n  \right)
- \Pr\left( A\nprec C(t_A)\right) \right| \leq b_1 \ e^{-b_2 n}\, .
$$ 
\end{lemma}

\bpr
We first observe that for any pattern $A$ and any positive integer
$\Delta$ such that $n+\Delta<t_A$, we have
$$ 
\Pr\left(A\nprec C(t_A)\backslash C_{n+\Delta}) \right)-
\Pr\left(A\nprec C(t_A)\right)= 
$$
$$
\Pr\left( A\prec C_{n+\Delta}\right)\leq (n+\Delta)^{d}\ \Pr(A)
$$
where we used Lemma \ref{volume_estim} to get the inequality.
For the sake of convenience, ``A is good'' stands for
$\forall {\bf x}\in\mathbb{Z}^d$ such that $0<|{\bf x}|<n/2$, we have
$(\tau_{\bf x}\sigma)_{C_n}\neq A$ for every $\sigma\in\C(A)$.
Now we use that $A$ is good to obtain
$$
\Pr\left(A\nprec C(t_A)\backslash C_{n+\Delta}, A\,\textup{is good}\
  \big| \ A\prec C_{n} \right)-
\Pr\left(A\nprec C(t_A)\backslash C_{n}, A\,\textup{is good}\ \big| \
  A\prec C_{n} \right)=
$$
\begin{equation}\label{P}
\Pr\left( \exists {\bf x}, n/2 < |{\bf x}| < n+\Delta\ : \
  \sigma_{C_n+{\bf x}\backslash C_n} = P_n^{{\bf x}}\ \big| \
  A\prec C_{n} \right)
\end{equation}
where $P_n^{{\bf x}}$ is a {\em fixed} pattern depending only on
$A=A_n$ and $\bf x$. Using the Gibbs property we obtain
$$
\aligned
\eqref{P} & \leq 
(n+\Delta)^d \ \sup_{|{\bf x}|>n/2}\ \sup_{\eta}
\Pr\left( \sigma_{C_n+{\bf x}\backslash C_n}=P_n^{{\bf x}}\ \big| \
  \eta_{C_{n}}=A \right) \}\\
& \leq  (n+\Delta)^d \ \sup_{|{\bf x}|>n/2}
\exp(-c\ |C_n+{\bf x}\backslash C_n|)\\
& \leq (n+\Delta)^d \ \exp(-c' n^d)
\endaligned
$$
where $c,c'$ are positive constants.
We now use the mixing property \eqref{mixfun} to get,
for any good pattern $A$~:
$$
\left|
\Pr\left(A\nprec C(t_A)\backslash C_{n+\Delta}\ \big| \ A\prec C_{n}
\right)
-
\Pr\left(A\nprec C(t_A)\backslash C_{n+\Delta} \right)
\right|
\leq
|C(t_A)|\ \phi(\Delta)\,. 
$$
Putting together the above estimates, with the choice $\Delta=n^{d+1}$
and using \eqref{exponmix}, yields
$$
\left|
\Pr\left( A\nprec C(t_A)\backslash C_n\ \big| \ A\prec C_n  \right)
- \Pr\left( A\nprec C(t_A)\right) \right| \leq 
$$
$$
(n+n^{d+1})^d \ e ^{-c'\ n^d} + C_1 \ e ^{c" n^d} e ^{-C_2 n^{d+1}}+
(n+n^{d+1})^d \ \Pr(A)
\, .
$$ 
This gives the desired result.
\epr

We also need the following lemma.

\begin{lemma}[Iteration Lemma for pattern repetitions]
\label{iterationlemmaforrepetition}
Let $t_A$ be such that $t_A^d\sim \Pr(A)^{-\vartheta}$, where
$\vartheta\in (0,1)$. For $i=2,\ldots k$,
let
$C_i (t_A)$ denote any collection of $k$ disjoints
cubes of the form ${\mathbf x}_i+ C(t_A)$. 
Assume also that $C_1 (t_A) ={\mathbf x}_1+ C(t_A)\backslash \{0\}$ is 
disjoint from $C_i (t_A)$, $i=2,\dots,k$.
Then 
we have the following inequality for all $k$:
\beq
&&
\Big| 
\Pr\left( A\nprec\bigcup_{i=1}^k C_i (t_A)  \ \vert \ A\prec C_n \right)
-
\Pr\left( A\nprec C(t_A)\right)^{k} \Big| \nonumber\\
&\leq& 
C_1 \exp\{- C_2 n \}  \
\left( 
\Pr\left( A\nprec C(t_A)\right)+C_1 \exp\{- C_2 n \}
\right)^{k}  . \nonumber
\eeq
\end{lemma}

\bpr 
Proceeding as in the proof of Lemma \ref{iterationlemma} we have:

\beq
\Big|\Pr\left( (A\prec C_n)\cap (A\nprec\bigcup_{i=1}^k C_i (t_A))
\right)
-
\Pr((A\prec C_n) \cap A\nprec C(t_A)\backslash C_n)\
\Pr\left( A\nprec C(t_A)\right)^{k-1} \Big| \nonumber
\leq 
\eeq
\beq
\Pr(A)^{1-\vartheta} \
\left( 
\Pr\left( A\nprec C(t_A)\right)+\Pr(A)^{1-\vartheta}
\right)^{k}  . \nonumber
\eeq
On the other hand, Lemma \ref{shortret} tells us that
\beq
\Big|
\Pr( A\nprec C(t_A)\backslash C_n \vert  A\prec C_n )
-
\Pr( A\nprec C(t_A) ) \Big| \nonumber
\le
b_1\ e^{-b_2 n}\, .
\eeq
\epr

The proof of \eqref{exponentiallawforrepetitions} in 
Theorem \ref{corollary-repetition} is now the same as
that of Theorem \ref{thm1}.
It remains to prove \eqref{probaofbadpatterns}:

\begin{lemma}[Probability of badly self-repeating
patterns]\label{bsr-lemma}
There exist $c,C>0$ such that
\be\label{bsr}
\Pr (\bee_n)\leq B e^{-bn^d}
\ee
\end{lemma}
\bpr
Put $C^+_n= C_n \cap (C_n+{\bf x})$ and $C^+_n= C_n \cap (C_n- {\bf x})$ .
By definition of $\bee_n$, we have the inequality:
\be
\Pr (\bee_n) \leq \Pr \left( \exists {\bf x}: |{\bf x}| \leq n/2 :
\sigma_{C^+_n ({\bf x})} =\sigma_{C^-_n ({\bf x})}\right).
\ee
Define the event $E_{{\bf x}}=\{\sigma:\sigma_{C^+_n ({\bf x})}
=\sigma_{C^-_n ({\bf x})}\}$.
If $\si\in\ E_{{\bf x}}$, then there exists disjoint sets $S^+_n ({\bf
x})$ and
$S^-_n ({\bf x})$ such that $\si_{S^+_n ({\bf x})}=\si_{S^-_n ({\bf x})}$
and
$|S^+_n ({\bf x})|,|S^-_n ({\bf x})| >\delta n^d$ for some positive
$\delta$.
Therefore, we have
\beq
\Pr (E_{\bf x}) &\leq &
\Pr \left(\si_{S^+_n ({\bf x})}=\si_{S^-_n ({\bf x})}\right)
\nonumber\\
&\leq &
\sup \left\{ \Pr \left(\si_{S^+_n ({\bf x})}=\eta|\si_{(S^+_n ({\bf
x}))^c}=\xi\right)
:\eta\in\Omega_{S^+_n ({\bf x})},\xi\in\Omega_{(S^+_n ({\bf
x}))^c}\right\}
\nonumber\\
&\leq &
\exp(-c' n^d)
\eeq
where in the last inequality we used the Gibbs property \eqref{gibbsineq}.
Finally, 
\be
\Pr (\bee_n) \leq \sum_{{\bf x}: |{\bf x}|< n/2} \Pr (E_{{\bf x}})
\leq B e^{-bn^d}.
\ee
\epr

\subsection{Proof of Theorem
\protect{\ref{corollary-OW}}}\label{proof-corollary-OW}

We start by showing the following summable upper-bound to
$$
\Pr\{\sigma~:\log(\r_n(\sigma)^d \Pr({\mathcal C}(\sigma_{C_n})))\geq \log
t\}\le
$$
$$
\sum_{\C_n\in \bee_n^c}\Pr(\C_n)\
\Pr\{\sigma~:\log(\r_n(\sigma)^d \Pr(\C_n))\geq \log t\ \vert\ \C_n\}
+\sum_{\C_n\in \bee_n}\Pr(\C_n)\,.
$$
From Theorem \ref{corollary-repetition} and Lemma \ref{bsr-lemma}
we get for all $t>0$
$$
\Pr\{\sigma~:\log(\r_n(\sigma)^d \Pr({\mathcal C}(\sigma_{C_n})))\geq \log
t \}\leq 
 (C' e^{-c'n^d}+ e^{-\Lambda_1 t}) +C e^{-cn^d}\, .
$$
Take $t=t_n=\log(n^{\epsilon})$, $\epsilon>\Lambda_1^{-1}$, to get
$$
\Pr\{\sigma~:\log(\r_n(\sigma)^d \Pr({\mathcal C}(\sigma_{C_n})))\geq
\log\log(n^{\epsilon})\}
\leq  C'e^{-c'n^d} + \frac1{n^{\epsilon\Lambda_1}}+ C e^{-cn^d}\,.
$$
An application of the Borel-Cantelli lemma leads to
$$
\log\left[(\r_n(\sigma))^d \Pr({\mathcal C}(\sigma_{C_n}))\right]\leq
\log\log(n^{\epsilon})\quad\textup{eventually a.s.}\,.
$$
For the lower bound first observe that Theorem \ref{corollary-repetition}
gives, for all $t>0$
$$
\Pr\{\sigma~:\log(\r_n(\sigma)^d \Pr({\mathcal C}(\sigma_{C_n})))\leq \log
t\}
\leq C' e^{-c'n^d}+(1-\exp(-\Lambda_2 t))+C e^{-cn^d}\, .
$$
Choose $t=t_n=n^{-\epsilon}$, $\epsilon>1$, to get, proceeding as before,
$$
\log\left[(\r_n(\sigma))^d \Pr({\mathcal C}(\sigma_{C_n}))\right]\geq
-\epsilon\log n\quad\textup{eventually a.s.}\,.
$$
Finally, let $\epsilon_0 = \max( \Lambda_1^{-1}, 1)$.

\subsection{Proof of Theorem \protect{\ref{waiting-time}}}\label{wtproof}

We first show that the strong approximation formula
\eqref{strong-approximation} holds
with $\w_n$ in place of $\r_n$ with respect to the measure
$\mathbb{Q}\times\Pr$.
We have the following identity:
$$
\int d\mathbb{Q}(\xi)\
\Pr\left\{\sigma: \t_{\xi_{C_n}}(\sigma)> \left(\frac{t}{\Pr({\mathcal
C}(\xi_{C_n}))}\right)^{1/d}\right\}=
$$
$$
(\mathbb{Q}\times \Pr)\left\{(\xi,\sigma): \w_n(\xi,\sigma)>
\left(\frac{t}{\Pr({\mathcal C}(\xi_{C_n}))}\right)^{1/d}\right\}
$$
This shows immediately that Theorem \ref{thm1}
is valid with $\w_n(\xi,\sigma)$ in place of $\t_{\sigma_{C_n}}(\xi)$
and $\mathbb{Q}\times\Pr$ in place of $\Pr$, hence so is Theorem
\ref{corollary-OW}.
Therefore for $\epsilon$ large enough, we obtain
\be\label{SAW}
-\epsilon\log n \leq\log\left[(\w_n(\xi,\sigma))^d \Pr({\mathcal
C}(\xi_{C_n}))\right]\leq \log\log n^\epsilon
\ee
for $\mathbb{Q}\times\Pr$-eventually almost every $(\xi,\sigma)$.
Write
$$
\log\left[(\w_n(\xi,\sigma))^d \Pr({\mathcal C}(\sigma_{C_n}))\right]=
d\log\w_n(\xi,\sigma)+\log\mathbb{Q}({\mathcal C}(\xi_{C_n}))-
\log\frac{\mathbb{Q}({\mathcal C}(\xi_{C_n}))}{\Pr({\mathcal
C}(\xi_{C_n}))}
$$
and use (\ref{SAW}). After division by $n^d$, we obtain (\ref{ASW})
since $\limn\frac1{n^d}\log\mathbb{Q}({\mathcal
C}(\sigma_{C_n}))=-s(\mathbb{Q})$, $\mathbb{Q}$-a.s. by the
Shannon-Mc Millan-Breiman theorem and
$\limn\frac1{n^d}\log\frac{\mathbb{Q}({\mathcal
C}(\xi_{C_n}))}{\Pr({\mathcal C}(\xi_{C_n}))}=s(\mathbb{Q}|\Pr)$,
$\mathbb{Q}$-a.s. (Proposition \ref{asre} in Section \ref{Gibbs}). 

\subsection{Proof of Theorem \protect{\ref{corollary-clt}} and
Theorem \protect{\ref{corollary-clt-waiting}}}\label{cltproof}

We use the strong approximation formula \eqref{strong-approximation} from
Theorem \ref{corollary-OW} 
to get
\begin{equation}
\frac{d\log\r_n(\sigma)+\log\Pr_\beta({\mathcal C}(\sigma_{C_n}))}
{n^{\frac{d}{2}}}\to 0\quad\textup{when}\,n\to\infty,\;\textup{for}\
\Pr_\beta-\textup{almost all}\, \sigma\, .
\end{equation}
Therefore, it suffices to see that in the high-temperature regime 
we have a central limit theorem
for $\{-\frac1{n^d}\log\Pr_\beta({\mathcal C}(\sigma_{C_n}))\}$. By a
standard
argument presented below (\ref{standard-fact}), one has
\be
\limn\frac1{n^d}\log\int \Pr_\beta ({\mathcal C}(\xi_{C_n}))^{-q}\
d\Pr(\xi)=
P((1-q)\beta U)+ (q-1)P(\beta U)\, ,
\ee
for all $q\in [0,\infty)$.
There exists $\beta_1>0$ such that 
for $|z|\leq \beta_1$ the maps $z\mapsto P(z U)$ and
\[
\Psi:z\mapsto \limn\frac1{n^d}\log\int \Pr_\beta ({\mathcal C}(\xi_{C_n}))^{-z}\
d\Pr(\xi)
\]
are analytic
see e.g. \cite{simon}, and \cite{dobshlos}. Therefore, if
$|(q-1)|\beta\leq \beta_1$, the map $q\mapsto P((1-q)\beta U)+ (q-1)
P(\beta U)$ is analytic, and equality holds for all $q\in\Comp$.

By Bryc's theorem \cite{bryc}, this implies the CLT for
$\{-\frac1{n^d}\log\Pr_\beta({\mathcal C}(\sigma_{C_n}))\}$
with variance $\theta^2$ given by
\be
\theta^2= \frac{d^2}{dq^2} \left(P((1-q)\beta U)\right)\big|_{q=0}
\ee
which is strictly positive by strict convexity of the pressure in the
analyticity regime.
The proof of Theorem \ref{corollary-clt-waiting} is the same once we
observe that
\begin{equation}
\frac{d\log\w_n(\sigma)+\log\Pr_\beta({\mathcal C}(\xi_{C_n}))}
{n^{\frac{d}{2}}}\to 0\quad\textup{when}\,n\to\infty,\;\textup{for}\
\Pr_\beta\times\Pr_\beta-\textup{almost all}\,
(\xi,\sigma)
\end{equation}
by using \eqref{SAW}.

\subsection{Proof of Theorem \protect{\ref{LD-waiting}}}\label{LD-proof}
Recall that for any Gibbs measure
$$
-\log\Pr(\sigma_{C_n})\sim \sum_{i\in C_n}\tau_i f_U(\sigma) +\log
Z_{C_n}\,
$$
and hence we have the identity
\be\label{standard-fact}
\limn \frac{1}{n^d} \log \sum_{\mathcal C_n} \Pr(\mathcal C_n)^{1-q}
= P( (1-q)U) -(1-q)P(U).
\ee
In the sequel, we are going to show that 
\begin{equation}\label{q_large}
   {\int \w_n^{qd} d\Pr\times\Pr}\approx \sum_{\mathcal C_n} \Pr( \mathcal
C_n)^{1-q},
\end{equation}
for $q>-1$, and
\begin{equation}\label{q_small}{\int \w_n^{qd} d\Pr\times\Pr}\approx
\sum_{\mathcal C_n} \Pr( \mathcal C_n)^{2},
\end{equation}
for $q\le -1$. Here $a_n\approx b_n$  means that
$\max\{a_n/b_n,b_n/a_n\}$ is bounded from above.
Clearly (\ref{q_large}) and (\ref{q_small}) imply (\ref{w-cases}).

Let $q> 0$. Then
\begin{align}\label{int-repres}
\int \w_n^{qd} d\Pr\times \Pr & =
\sum_{\mathcal C_n} \Pr(\mathcal C_n) \int 
\t_{\mathcal C_n}^{qd}(\sigma) d\Pr(\sigma) \\
& = q\sum_{\mathcal C_n} \Pr(\mathcal C_n)^{1-q}
\int_{\Pr(\mathcal C_n)}^{\infty} t^{q-1} \Pr \left\{ \t_{\mathcal
C_n}^d\ge \frac 
t{\Pr(\mathcal C_n)}\right\} dt.
\end{align} 
By Theorem\ref{thm1},
there exist positive constants $A,B$ such
that for any $t>0$ one has
$$
\Pr \left\{ \t_{\mathcal C_n}^d\ge \frac 
t{\Pr(\mathcal C_n)}\right\}\le A e^{-Bt}.
$$
Theorem\ref{thm1} also easily gives the lower bound~:
$$
\int_{\Pr(\mathcal{C}_n)}^{\infty} t^{q-1} \Pr \left\{
\t_{\mathcal{C}_n}^d>
\frac t{\Pr(\mathcal{C}_n}\right\} dt\geq
K' - C\exp(-c n)\ K^{\prime\prime}
$$
where
$0<K':=\int_1^\infty t^{q-1}\ e^{-\Lambda_2 t}\ dt<\infty$ and
$0<K^{\prime\prime}:=\int_0^\infty t^{q-1}\ e^{-\Lambda_1 t}\ dt<\infty$. 
For $n$ large enough, $K' - C\exp(-c n)K^{\prime\prime}$ is strictly positive.
Therefore we obtain
$$K_1\sum_{\mathcal C_n} \Pr(\mathcal C_n)^{1-q}\le \int \w_n^{qd}
d\Pr\times \Pr
\le K_2\sum_{\mathcal C_n} \Pr(\mathcal C_n)^{1-q},
$$
where $$
 K_1 := q \ (K' - C\exp(-c n_0)\ K^{\prime\prime})\, ,
 \quad
 K_2 := qA \int_{0}^{\infty} t^{q-1} e^{-Bt} dt.
$$
This establishes (\ref{q_large}) for $q\ge 0$. The case $q=0$ is
trivial.

Let now $q\in (-1,0)$.

  \begin{align}\label{int-repres1}
\int \w_n^{-|q|d}\ d\Pr\times \Pr & =
\sum_{\mathcal C_n} \Pr(\mathcal C_n) \int 
\t_{\mathcal C_n}^{-|q|d}(\sigma)\ d\Pr(\sigma) \\
& = \sum_{\mathcal C_n} \Pr(\mathcal C_n) 
\int_{0}^{1} \Pr \left\{ \t_{\mathcal C_n}^{-|q|d}\ge t\right\} dt\\
&= |q|\sum_{\mathcal C_n} \Pr(\mathcal C_n)^{1+|q|} 
\int_{\Pr(\mathcal C_n)}^{\infty} t^{-|q|-1}\ \Pr \left\{ \t_{\mathcal
C_n}^d\le \frac 
t{\Pr(\mathcal C_n)}\right\} dt.
\end{align} 

The last integral is bounded from above by the integral where
$\Pr(\mathcal{C}_n)$ replaced by $1$ in the integration domain. 
From Theorem \ref{thm1},  we get the following lower bound,
for every $t>0$:
$$
\Pr \left\{ \t_{\mathcal C_n}^d\le \frac  t{\Pr(\mathcal C_n)}
\right\}
\ge
1- e^{-\Lambda_{1}t}- C'\Pr({\mathcal C_n})^\rho e^{-C''t}
$$
The number 
$$
K_1':= |q| \int_{1}^{\infty} t^{-|q|-1} \left( 1- e^{-\Lambda_{1}t}-
  C'\Pr({\mathcal C_n})^\rho e^{-C''t}\right)\ dt
$$
is finite and strictly positive for $n$ large enough.

Now, putting $0$ instead of $\Pr(\mathcal{C}_n)$ gives an upper
bound to the integral upon condideration. We use 
Lemma \ref{pospar} to get immediately
$$
\Pr \left\{ \t_{\mathcal C_n}^d\le \frac t{\Pr(\mathcal C_n)} \right\}
\le 
1- e^{ -\Lambda_{2}t}\ .
$$
provided that $t\leq \frac1{2}$. We
have
$$
\int_{0}^{\infty} t^{-|q|-1}\ \Pr \left\{ \t_{\mathcal C_n}^d\le \frac 
t{\Pr(\mathcal C_n)}\right\} dt \leq
$$
$$
\int_{0}^{\frac1{2}} t^{-|q|-1}\ \Pr \left\{ \t_{\mathcal C_n}^d\le \frac 
t{\Pr(\mathcal C_n)}\right\} dt + 
\int_{\frac1{2}}^{\infty} t^{-|q|-1}\ dt \leq
$$
$$
\frac{\Lambda_2\ 2^{1-|q|}}{1-|q|}+ \frac{2^{-|q|}}{|q|}=: K_2'<\infty\ .
$$
Hence, we conclude that for
$n$ large enough
$$
K_1'\sum_{\mathcal C_n} \Pr(\mathcal C_n)^{1+|q|}\le \int \w_n^{-|q|d}
d\Pr\times \Pr
\le |q|\ K_2'\sum_{\mathcal C_n} \Pr(\mathcal C_n)^{1+|q|}\ .
$$
Therefore we obtain (\ref{q_large}) for $q\in(-1,0)$.

Finally, let us consider the remaining case $q\le-1$. Then
for sufficiently large $n$ (such that $\Pr(\mathcal{C}_n)<1/2$) one has
\begin{align*}
\int \w_n^{-|q|d} d\Pr\times \Pr&= |q|\sum_{\mathcal C_n} \Pr(\mathcal
C_n)^{1+|q|}\ 
\int_{\Pr(\mathcal C_n)}^{\infty} t^{-|q|-1}\ \Pr \left\{ \t_{\mathcal
C_n}^d\le \frac 
t{\Pr(\mathcal C_n)}\right\} dt\\
&=|q|\sum_{\mathcal C_n} \Pr(\mathcal C_n)^{1+|q|}\ 
\Bigl[
\int_{\Pr(\mathcal C_n)}^{\frac1{2}} 
+\int_{\frac1{2}}^{\infty} \Bigr]\
t^{-|q|-1} \Pr \left\{ \t_{\mathcal C_n}^d\le \frac 
t{\Pr(\mathcal C_n)}\right\} dt\\
&=|q|\sum_{\mathcal C_n} \Pr(\mathcal C_n)^{1+|q|} \
\left[ \; I_1(n,\mathcal C_n)  + I_2(n,\mathcal C_n)\, \right]. 
\end{align*}
Clearly the second integral $I_2(n,\mathcal C_n)$ is uniformly bounded in
$n$. Indeed, 
$$
I_2(n,\mathcal C_n)\le \int_{\frac1{2}}^\infty \frac 1{t^{1+|q|}}\ dt
<+\infty.
$$
However, the first integral $I_1(n,\mathcal C_n)$ is diverging in the
limit $n\to\infty$. 
Therefore the limiting behavior as $n\to\infty$ is determined by
$$
|q|\sum_{\mathcal C_n} \Pr(\mathcal C_n)^{1+|q|} 
I_1(n,\mathcal C_n)=
|q|\sum_{\mathcal C_n} \Pr(\mathcal C_n)^{1+|q|} 
\int_{\Pr(\mathcal C_n)}^{\frac1{2}} 
t^{-1-|q|} \Pr \left\{ \t_{\mathcal C_n}^d\le \frac 
t{\Pr(\mathcal C_n)}\right\}\ dt.
$$
We again use Lemma \ref{pospar} to get
$$
1- e^{ -\Lambda_{1}t}\le
\Pr \left\{ \t_{\mathcal C_n}^d\le \frac t{\Pr(\mathcal C_n)} \right\}
\le 
1- e^{ -\Lambda_{2}t}\ .
$$
provided that $t\leq \frac1{2}$.
Hence, using the Gibbs property \eqref{gibbsineq}, we have
$$
I_1(n,\mathcal{C}_n)\leq \Lambda_2 \int_{\Pr(\mathcal C_n)}^{\frac1{2}}
t^{-|q|}\ dt
\leq
\frac{\Lambda_2 (1-2^{|q|-1}C' e^{-c' n^d}) }{|q|-1}\
\Pr(\mathcal{C}_n)^{-|q|+1}
$$
where we used the fact that for all $\kappa\in\R$, $1-e^{-\kappa}\le
\kappa$. Notice that for $n$ large enough, the term between parentheses is
strictly positive.
Now, using the fact that $1-e^{-\kappa}\ge \kappa/2$ for any
$\kappa\in[0,1]$, and remembering that $\Lambda_{1} /2 \le 1$
(\footnote{Indeed, $\Lambda_1 \le \Lambda_2=2$, see the end of the proof
of Lemma \ref{pospar}.}), and using again the Gibbs property
\eqref{gibbsineq}, we obtain
$$
I_1(n,\mathcal{C}_n)\geq \frac{\Lambda_1 (1-2^{|q|-1} C e^{-c
n^d})}{2(|q|-1)}\ \Pr(\mathcal{C}_n)^{-|q|+1}
$$ 
where the term between parenthese is strictly positive provided that
$n$ is sufficiently large.
Therefore, for $n$ large enough, we end up with
$$
\frac{|q|\Lambda_1 (1\!-\!2^{|q|-1}C e^{-c n^d})}{2(|q|-1)}\
\!\!\sum_{\mathcal C_n} \Pr(\mathcal C_n)^{2} 
\le \int \w_n^{-|q|d} d\Pr\times\Pr \le
\frac{2|q|\Lambda_2 (1\!-\!2^{|q|-1}C' e^{-c' n^d})}{|q|-1}\
\!\!\sum_{\mathcal C_n} \Pr(\mathcal C_n)^{2}.
$$
(Notice that L'H{\^o}pital's rule shows that there is no problem at $q=-1$.)
Thus, we obtain (\ref{q_small}), which finishes the proof.


\end{document}